\documentclass[11pt,reqno,letterpaper]{amsart}
\usepackage{amsmath,amsthm,amssymb,mathtools,bbm}
\usepackage{graphicx}
\usepackage{hyperref}
\usepackage{breakurl}
\usepackage{bookmark}
\hypersetup{pdfstartview={FitH}}

\theoremstyle{plain}
\newtheorem{theorem}{Theorem}

\theoremstyle{definition}
\newtheorem{remark}[theorem]{Remark}

\numberwithin{equation}{section}

\newcommand{\R}{\mathbb{R}}
\newcommand{\Z}{\mathbb{Z}}
\newcommand{\1}{\mathbbm{1}}
\newcommand{\norm}[2]{\lVert #1\rVert_{\textup{L}^{#2}}}

\renewcommand{\le}{\leqslant}
\renewcommand{\ge}{\geqslant}

\title{Sharpness of the Martin--O'Bryant autoconvolution inequality}

\author[T. Crmari\'{c}]{Ton\'{c}i Crmari\'{c}}
\address[T.\,C.]{Department of Mathematics, Faculty of Science, University of Split, Ru\dj{}era Bo\v{s}kovi\'{c}a 33, 21000 Split, Croatia}
\email{tcrmaric@pmfst.hr}

\author[V. Kova\v{c}]{Vjekoslav Kova\v{c}}
\address[V.\,K.]{University of Zagreb Faculty of Science, Department of Mathematics, Bijeni\v{c}ka cesta 30, 10000 Zagreb, Croatia}
\email{vjekovac@math.hr}

\author[J. Rup\v{c}i\'{c}]{Jelena Rup\v{c}i\'{c}}
\address[J.\,R.]{Faculty of Transport and Traffic Sciences, University of Zagreb, Vukeli\'{c}eva 4, 10000 Zagreb, Croatia}
\email{jrupcic@fpz.unizg.hr}


\begin{document}

\begin{abstract}
We show that there does not exist a constant $c<1$ such that $\norm{f\ast f}{2}^2\le c\,\norm{f\ast f}{1}\norm{f\ast f}{\infty}$ holds for every nonnegative $f\in\textup{L}^1(\R)\cap\textup{L}^2(\R)$.
\end{abstract}

\maketitle


\section{Introduction}

For every nonzero nonnegative $f\in\textup{L}^1(\R)\cap\textup{L}^2(\R)$, set
\[ \mathfrak{R}(f) := \frac{\norm{f\ast f}{2}^2}{\norm{f\ast f}{1}\norm{f\ast f}{\infty}}. \]
Let $c$ be the supremum of the ratios $\mathfrak{R}(f)$ over all such $f$. Clearly, $c\le1$ holds by the straightforward inequality
\[ \norm{f\ast f}{2}^2 \le \norm{f\ast f}{1} \norm{f\ast f}{\infty}. \]
Martin and O'Bryant used this inequality in a proof \cite[\S4]{MO02} related to a continuous variant of a question from Ramsey theory, asked if it can be sharpened, and conjectured \cite[Conj.~4.4]{MO02} that the optimal upper bound on the ratios $\mathfrak{R}(f)$ was $(\log 16)/\pi$.
They restated the conjecture in \cite[Conj.~2.9]{MO07} and \cite[Conj.~5.2]{MO}, but Matolcsi and Vinuesa \cite{MV} disproved the proposed bound. The value of $c$ remained open and was subsequently discussed by O'Bryant \cite{MOquestion} and by Barnard and Steinerberger \cite[\S1.4]{BS}.

\begin{table}[htbp]
\centering
\caption{Selected lower bounds for $c$. Some numerical values are rounded.}
\label{tab:records}
\small
\begin{tabular}{@{}lll@{}}
\textbf{Date} & \textbf{Construction} & \textbf{Value} \\
Oct.\ 2002 & Martin and O'Bryant \cite{MO02} & $0.88254240$ \\
July 2009 & Matolcsi and Vinuesa, 19 steps \cite[arXiv v1]{MV} & $0.88379269$ \\
Sept.\ 2009 & Matolcsi and Vinuesa, 20 steps \cite[arXiv v2]{MV} & $0.88922617$ \\
May 2025 & Novikov et al. (AlphaEvolve), 50 steps \cite{AE} & $0.89620000$ \\
June 2025 & Boyer and Li, 575 steps \cite[arXiv v1]{BL} & $0.90156200$ \\
Aug.\ 2025 & Jaech and Joseph, 2399 steps \cite[\S3]{JJ} & $0.94136000$ \\
Nov.\ 2025 & Georgiev et al. (AlphaEvolve), 50000 steps \cite{GGTW} & $0.96100000$ \\
Feb.\ 2026 & Kravatskiy et al. (ImprovEvolve) \cite{Improv} & $0.96258000$ \\
Mar.\ 2026 & Kang (ClaudeExplorer) \cite{KangMarch} & $0.96272000$ \\
May 2026 & Kang (Stochastic Agent Descent) \cite{Kang} & $0.96273000$ \\
June 23, 2026 & Hyra \cite{Hyra,Arena,Russell} & $0.96290110$ \\
July 19, 2026 & ClaudeExplorer \cite{Arena} & $0.96343273$ \\
July 21, 2026 & MLAI-Yonsei \cite{Arena} & $0.96345572$ \\
July 22, 2026 & ClaudeExplorer \cite{Arena} & $0.96358811$ \\
November 14, 2026 & MLAI-Yonsei \cite{Arena} & $0.96366156$ \\
This paper & Theorem~\ref{thm:main} & $1$ \\
\end{tabular}
\end{table}

Table~\ref{tab:records} traces selected records on the lower bounds for $c$, from the original conjecture to the highest score currently displayed on EinsteinArena \cite{Arena} at the time of writing of this manuscript.
Also note that Finch \cite[\S4.21]{Finch} surveyed the early bounds.

The following theorem establishes that $c=1$.

\begin{theorem}\label{thm:main}
For every $\varepsilon>0$ there exist a nonnegative integer $d$ and nonnegative numbers $p_0,\ldots,p_d$, not all zero, such that the step function
\[ f = \sum_{j=0}^d p_j\1_{[j,j+1)} \]
satisfies $\mathfrak{R}(f)>1-\varepsilon$.
\end{theorem}

For real-valued functions allowed to change sign, the problem is easier: de Dios Pont and Madrid \cite[Thm.~5.1]{DM} obtained the sharp constant $1$ by taking a square root of the Fourier transform.
Chung, Du, and Wesley \cite[\S6.12]{CDW} proved that restricting $f$ to $\{0,1\}$-valued step functions on finite grids gives the same supremum $c$, a reduction suggested by Tao in \cite{MOquestion}.

Rogers and one of the present authors \cite[\S7.1]{KR26} studied a similar problem on whether the autocorrelation of a nonnegative function of integral $1$ can get close to an indicator function in the $\textup{L}^1$ norm. A relatively simple recursive construction sufficed in \cite{KR26}, but the present problem is more subtle in the sense that $f\ast f$ also needs to be controlled in the $\textup{L}^\infty$ norm.

We present a proof of Theorem~\ref{thm:main} in Section~\ref{sec:proof}. In Section~\ref{sec:family} we illustrate the construction and provide an explicit choice of parameters depending on a single positive integer $m$. This gives a sequence of step functions $(f_m)$ satisfying $\mathfrak{R}(f_m)\to 1$, which is entirely concrete but difficult to compute due to rapidly growing numbers of steps.

We also formally verified Theorem~\ref{thm:main} in the proof assistant Lean 4 using its library Mathlib. This formalization is listed in the Palomar registry \cite{palomar}.


\section{Proof of Theorem~\ref{thm:main}}
\label{sec:proof}
The basic idea is to construct $f$ so that $f\ast f$ is has ``plateaus'' above many disjoint intervals. Their heights will be close to the global maximum, and together they will carry almost all of the mass.

Just as in the proof of \cite[Cor.~1.5]{MO}, it will turn out that step functions are conveniently encoded by polynomials.
For a nonzero polynomial $P(z)=\sum_{j=0}^d p_jz^j$ with nonnegative coefficients $p_j$, set
\[ f_P(x) := \sum_{j=0}^d p_j\1_{[j,j+1)}(x). \]
Extend $p_j$ by zero to all $j\in\Z$. Then the coefficient of $z^j$ in $P(z)^2$ is
\[ a_j = [z^j]P(z)^2=\sum_{k\in\Z}p_kp_{j-k} , \]
with $a_j=0$ for $j<0$ and $j>2d$.
The autoconvolution $\1_{[0,1)}\ast\1_{[0,1)}$ is a tent of height $1$ supported on $[0,2]$, and
\[ (f_P\ast f_P)(x) = \sum_{j\in\Z} a_j (\1_{[0,1)}\ast\1_{[0,1)})(x-j). \]
Consequently $f_P\ast f_P$ is piecewise linear with $(f_P\ast f_P)(j+1)=a_j$.
It is also nonnegative, so
\begin{equation}\label{eq:P_L1}
\norm{f_P\ast f_P}{1} = \sum_{j\in\Z}a_j=P(1)^2
\end{equation}
and
\begin{equation}\label{eq:P_Linf}
\norm{f_P\ast f_P}{\infty} = \max_{j\in\Z}a_j .
\end{equation}
Moreover, if $a_j=\cdots=a_{j+s}=H$, then $f_P\ast f_P$ equals $H$ throughout the segment $[j+1,j+s+1]$.
This plateau contributes $sH$ to $\int_{\R}f_P\ast f_P$ and $sH^2$ to $\int_{\R}(f_P\ast f_P)^2$.

We next construct polynomials $P_i$ for which $\sum_{i=0}^s P_i(z) P_{s-i}(z)$ has all coefficients equal to $1$.
Define $P_i(z)$ by the formal power series in $t$ as
\begin{equation}\label{eq:generating_fn}
\sum_{i=0}^{\infty}P_i(z)t^i=\frac{1}{\sqrt{(1-t)(1-zt)}}.
\end{equation}
The binomial expansion $(1-t)^{-1/2}=\sum_{k=0}^\infty 4^{-k}\binom{2k}{k}t^k$ gives an explicit formula
\begin{equation}\label{eq:binomial}
P_i(z) = 4^{-i}\sum_{h=0}^i \binom{2h}{h}\binom{2(i-h)}{i-h}z^h .
\end{equation}
In particular, $P_i$ has degree $i$ and positive coefficients.
At $z=1$ the generating function \eqref{eq:generating_fn} is $(1-t)^{-1} = \sum_{i=0}^\infty t^i$, so
\begin{equation}\label{eq:P1is1}
P_i(1)=1.
\end{equation}
Squaring the defining series \eqref{eq:generating_fn} gives
\[ \sum_{s=0}^\infty t^s \sum_{i=0}^s P_i(z)P_{s-i}(z) = \frac{1}{(1-t)(1-zt)} = \sum_{a,b\ge0}z^b t^{a+b} = \sum_{s=0}^\infty t^s \sum_{h=0}^s z^h. \]
Hence
\begin{equation}\label{eq:runs}
\sum_{i=0}^s P_i(z)P_{s-i}(z) = \sum_{h=0}^s z^h.
\end{equation}

In the next step we separate the coefficient blocks by their digits in a certain base $D$.
Fix a rational number $q\in(0,1)$ and positive integers $K<n$. Choose integers $r>b\ge 8n^3$ such that $b/r=q$, let $D=4n+1$, and set
\[ U(z) := \sum_{j=1}^{r}z^{D^{2j-1}},\quad
V(z) := \sum_{j=1}^{b}z^{D^{2j}},\quad
F(z) := \sum_{i=0}^{K}U(z)^{n-i}V(z)^{n+i}P_i(z). \]
We will use the step function $f_F$ to prove the theorem.
Squaring $F(z)$ gives
\[ F(z)^2 = \sum_{s=0}^{2K} W_s(z) B_s(z), \]
where
\begin{align}
W_s(z) & := U(z)^{2n-s} V(z)^{2n+s}, \label{eq:Wsfactors} \\
B_s(z) & := \sum_{i=\max\{0,s-K\}}^{\min\{K,s\}} P_i(z)P_{s-i}(z). \nonumber
\end{align}
For $0\le s\le K$, formula \eqref{eq:runs} gives
\begin{equation}\label{eq:property_B}
B_s(z)=1+z+\cdots+z^s,
\end{equation}
thanks to our choice of the polynomials $P_i$. For $K<s\le2K$, the index set is a nonempty subset of $\{0,\ldots,s\}$.
Each summand is a polynomial in $z$ of degree $s$ with a positive coefficient at every degree from $0$ through $s$.
Since the full sum in \eqref{eq:runs} has every such coefficient equal to $1$, the coefficient of $z^h$ in $B_s(z)$ satisfies
\[ 0<[z^h]B_s(z) \le1 \]
for $0\le s\le2K$ and $0\le h\le s$.

In order to investigate a typical term in the expansion of \eqref{eq:Wsfactors}, consider a choice of one monomial in each of the $2n-s$ factors $U(z)$ and
$2n+s$ factors $V(z)$.
If $\alpha_1,\ldots,\alpha_r$ and $\beta_1,\ldots,\beta_b$ are the multiplicities of the monomials from the two groups, its exponent is
\[ k := \sum_{j=1}^{r}\alpha_jD^{2j-1}+\sum_{j=1}^{b}\beta_jD^{2j},\quad
\sum_{j=1}^{r}\alpha_j=2n-s,\quad
\sum_{j=1}^{b}\beta_j=2n+s . \]
Every multiplicity, $\alpha_j$ or $\beta_j$, is at most $2n+2K<4n<D$. These are therefore the actual digits of $k$ in base $D$, with the units digit being zero. In particular, $k$ determines both multiplicity tuples, $(\alpha_1,\ldots,\alpha_r)$ and $(\beta_1,\ldots,\beta_b)$, and also determines $s$, since the sum of its digits in positions $2,4,\ldots,2b$ is precisely $2n+s$.
For $0\le h\le s\le 2K<D$, the number $k+h$ has the same higher digits as $k$ and the units digit equal to $h$. Thus, the number $k+h$ uniquely determines $k$, $h$, and $s$.
In the square $F(z)^2$ this gives
\[ [z^{k+h}]F(z)^2=([z^k]W_s(z))([z^h]B_s(z)). \]
As $h$ runs from $0$ to $s$ these coefficients form one ``block.'' Distinct blocks start with distinct multiples of $D$, at distance at least $D$ apart. For $s\le K$ all the coefficients of the block equal $[z^k]W_s(z)$, thanks to the observed property of \eqref{eq:property_B}, giving a plateau on $[k+1,k+s+1]$. These plateau intervals are pairwise disjoint.

Finally we now count the blocks giving maximal plateau height at each level, i.e., the exponents at which $W_s$ attains its largest coefficient.
The number of ordered choices giving $k$, and hence its coefficient, is
\[ [z^k]W_s(z) = \frac{(2n-s)!}{\prod_{j=1}^{r}\alpha_j!} \,\frac{(2n+s)!}{\prod_{j=1}^{b}\beta_j!}. \]
Namely, the two multinomial coefficients above count the arrangements of the chosen monomials in their respective groups of factors. The whole product is at most
\[ A_s := (2n-s)!\,(2n+s)! , \]
with equality precisely when all the multiplicities $\alpha_j,\beta_j$ are $0$ or $1$. (Since $r,b \ge 8n^3 \ge 4n$, such choices exist.) If $N_s$ denotes the number of exponents $k$ giving equality, choosing the occupied digit positions gives
\[ N_s = \binom{r}{2n-s}\binom{b}{2n+s}. \]
Each of these $N_s$ exponents $k$ accounts for $A_s$ ordered choices. The total number of ordered choices is, on the other hand,
\[ W_s(1) = U(1)^{2n-s} V(1)^{2n+s} = r^{2n-s}b^{2n+s}=(rb)^{2n}q^s. \]

To estimate $N_s A_s/W_s(1)$, we will use an elementary inequality
\begin{equation}\label{eq:elem_ineq}
\prod_{j=1}^N (1-x_j) \ge 1 - \sum_{j=1}^N x_j,
\end{equation}
which holds for $x_1,\ldots,x_N\in[0,1]$.
Recall $s\le 2K<2n$ and write
{\allowdisplaybreaks\begin{align*}
\frac{N_s A_s}{W_s(1)} & = \frac{r(r-1)\cdots(r-2n+s+1) b(b-1)\cdots(b-2n-s+1)}{r^{2n-s} b^{2n+s}} \\
& = \biggl(\prod_{j=1}^{2n-s-1} \Bigl(1-\frac{j}{r}\Bigr)\biggr) \biggl(\prod_{j=1}^{2n+s-1} \Bigl(1-\frac{j}{b}\Bigr)\biggr)
\stackrel{\eqref{eq:elem_ineq}}{\ge} 1 - \sum_{j=1}^{2n-s-1} \frac{j}{r} - \sum_{j=1}^{2n+s-1} \frac{j}{b} \\
& = 1 - \frac{(2n-s-1)(2n-s)}{2r} - \frac{(2n+s-1)(2n+s)}{2b}
\ge 1 - \frac{(4n-1)4n}{16n^3} > 1 - \frac{1}{n}.
\end{align*}}
Therefore,
\begin{equation}\label{eq:mass}
N_s A_s \ge \Bigl(1-\frac{1}{n}\Bigr) W_s(1) = \Bigl(1-\frac{1}{n}\Bigr)(rb)^{2n}q^s.
\end{equation}
We can say that repeated digit choices account for at most a fraction $1/n$ of all choices, so they become negligible for large $n$.

It remains to sum the plateau contributions.
For $0\le s<2K$ we have $A_{s+1}>A_s$, so that $A_0\le A_s\le A_{2K}$. Every coefficient of $F(z)^2$ contributed by a block at the $s$th level is at most $A_s$, so
\[ \norm{f_F\ast f_F}{\infty} \stackrel{\eqref{eq:P_Linf}}{\le} A_{2K}. \]
Its integral is
\[ \norm{f_F\ast f_F}{1} \stackrel{\eqref{eq:P_L1}}{=} F(1)^2,\quad F(1)=\sum_{i=0}^K r^{n-i}b^{n+i} P_i(1) \stackrel{\eqref{eq:P1is1}}{=} (rb)^n \sum_{i=0}^Kq^i. \]
For $0\le s\le K$, each of the $N_s$ blocks counted in \eqref{eq:mass} supplies a plateau of height $A_s$ and length $s$.
Let $E$ be the union of these plateau intervals. Disjointness and \eqref{eq:mass} give
\[ \int_E f_F\ast f_F=\sum_{s=0}^K sN_s A_s \ge \Bigl(1-\frac{1}{n}\Bigr)(rb)^{2n}\sum_{s=0}^K s q^s. \]
Since the height of $f_F\ast f_F$ on each of these intervals is at least $A_0$, we have
\[ \int_{\R} (f_F\ast f_F)^2 \ge \int_E (f_F\ast f_F)^2 \ge A_0 \int_E f_F\ast f_F, \]
so that
\[ \norm{f_F\ast f_F}{2}^2 \ge A_0 \Bigl(1-\frac{1}{n}\Bigr)(rb)^{2n}\sum_{s=0}^K s q^s . \]
Combining this estimate with
\[ \norm{f_F\ast f_F}{1}\norm{f_F\ast f_F}{\infty} \le A_{2K} (rb)^{2n} \Bigl(\sum_{i=0}^K q^i\Bigr)^2 \]
and canceling $(rb)^{2n}$ yields
\begin{equation}\label{eq:lower}
\mathfrak{R}(f_F)\ge\Bigl(1-\frac{1}{n}\Bigr)\frac{A_0}{A_{2K}}\,
 \frac{\sum_{s=0}^K s q^s}{\bigl(\sum_{i=0}^K q^i\bigr)^2}.
\end{equation}
For a fixed $K$,
\[ \frac{A_0}{A_{2K}} = \frac{(2n)!^2}{(2n-2K)!(2n+2K)!}
 = \prod_{j=1}^{2K}\frac{2n-2K+j}{2n+j} \longrightarrow 1 \]
as $n\to\infty$. In words, the ratio between the guaranteed plateau height $A_0$ and the global upper bound $A_{2K}$ approaches $1$.
Thus the right-hand side of \eqref{eq:lower} tends to
\[ \frac{\sum_{s=0}^K s q^s}{\bigl(\sum_{i=0}^K q^i\bigr)^2}. \]
As $K\to\infty$, the numerator tends to $q/(1-q)^2$ and the sum in the denominator tends to $1/(1-q)$, so this quotient tends to $q$.
By choosing a rational number $q<1$ sufficiently close to $1$, then taking a large $K$, and finally choosing a large $n>K$, we can get the right-hand side of \eqref{eq:lower} arbitrarily close to $1$.
This proves Theorem~\ref{thm:main}.
\hfill$\square$

\begin{remark}
It is interesting to remark that the number $0<q<1$ served as a ``bias'' between shifts by odd and even powers of $D$. We let $q\to 1$ only at the very end of the proof. If we immediately took $q=1$, i.e., $b=r$, then there would be no guarantee that the numbers $\mathfrak{R}(f_F)$ get close to $1$. Namely, the last factor on the right-hand side of \eqref{eq:lower} would equal $K/(2(K+1))$, so taking $n\to\infty$ with $K$ fixed and then $K\to\infty$ would therefore give $1/2$ for this lower bound.
\end{remark}

\begin{remark}
One might wonder if the polynomials $P_i$ used in the proof belong to any special classical family. If we denote the Legendre polynomials by $\widetilde{P}_i$, then their generating function \cite[Eq.~14.7.19]{DLMF} reads
\[ \sum_{i=0}^\infty \widetilde{P}_i(x) h^i = \frac{1}{\sqrt{1-2xh+h^2}} \]
for $x,h\in(-1,1)$. For each $i$, the coefficients of $h^i$ on both sides are polynomials in $x$. Their equality on $(-1,1)$ therefore gives equality for every real $x$.
Substituting $x=(z^{1/2}+z^{-1/2})/2$ and $h=z^{1/2}t$ for $z,t\in(0,1)$ and comparing the coefficients to those in \eqref{eq:generating_fn}, we obtain
\[ P_i(z) = z^{i/2} \widetilde{P}_i\Bigl(\frac{z^{1/2}+z^{-1/2}}{2}\Bigr). \]
Thus, the polynomials $P_i$ are in fact obtained from the Legendre polynomials by a nonlinear change of variables.
\end{remark}


\section{An explicit sequence of examples}
\label{sec:family}

For any integer $m\ge2$, take
\[ q=1-\frac1m,\quad K=2m^3,\quad n=8m^8,\quad r=8mn^3,\quad b=8(m-1)n^3 . \]
Write $f_m=f_F$ for the resulting step function. For the factorial factor in \eqref{eq:lower}, we have
\[ \frac{A_0}{A_{2K}} = \prod_{j=1}^{2K} \Bigl(1-\frac{2K}{2n+j}\Bigr)
\stackrel{\eqref{eq:elem_ineq}}{\ge} 1-\sum_{j=1}^{2K}\frac{2K}{2n+j} \ge 1-\frac{2K^2}{n}=1-\frac{1}{m^2}. \]
Substituting this and
{\allowdisplaybreaks\begin{align*}
\frac{\sum_{s=0}^K s q^s}{\bigl(\sum_{i=0}^K q^i\bigr)^2}
& \ge \frac{q/(1-q)^2 - \sum_{s=K+1}^\infty s q^s}{1/(1-q)^2} \\
& \ge q - \frac{(1-q)^2}{K} \sum_{s=0}^\infty s (s-1) q^s \\
& = q - \frac{2q^2}{K(1-q)} = q - \frac{q^2}{m^2} \ge q \Bigl(1-\frac{1}{m^2}\Bigr)
\end{align*}}
into \eqref{eq:lower} yields
\begin{equation}\label{eq:step_est}
\mathfrak{R}(f_m) \ge \Bigl(1-\frac{1}{8m^8}\Bigr) \Bigl(1-\frac{1}{m}\Bigr) \Bigl(1-\frac{1}{m^2}\Bigr)^2 .
\end{equation}
In particular, we again see that $\lim_{m\to\infty}\mathfrak{R}(f_m)=1$.

\begin{figure}[htbp]
\begin{center}
\includegraphics[width=0.9\linewidth]{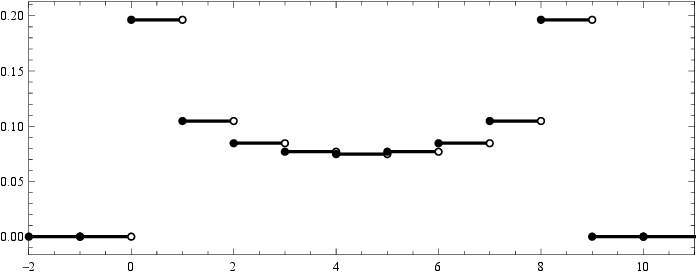}
\end{center}
\caption{Nine consecutive positive steps of $t\mapsto f_2(J_8+t)$. Created in Mathematica \cite{Mathematica}.}
\label{fig:m2-step}
\end{figure}

For \underline{$m=2$}, write $p_j=[z^j]F(z)$. The parameters are $K=16$, $n=2048$, and $D=8193$. For
$0\le i\le K$ let
\[ J_i := (n-i)D+(n+i)D^2 . \]
At $J_i+h$ with $0\le h\le i$, all $n-i$ choices from $U(z)$ use the first odd position, and all $n+i$ choices from $V(z)$ use the first even position. Both multinomial coefficients are $1$. For $i=8$, formula \eqref{eq:binomial} gives $p_{J_8+h}=[z^h]P_8(z)$ for $h=0,\ldots,8$ as the tuple
\[ (p_{J_8+h})_{h=0}^{8}
= \frac{1}{32768}(6435,3432,2772,2520,2450,2520,2772,3432,6435). \]
The cells immediately before and after the block have height zero. Figure~\ref{fig:m2-step} shows the graph of $f_2$ relative to the left endpoint of the chosen block, $J_8$.

\begin{figure}[htbp]
\begin{center}
\includegraphics[width=0.9\linewidth]{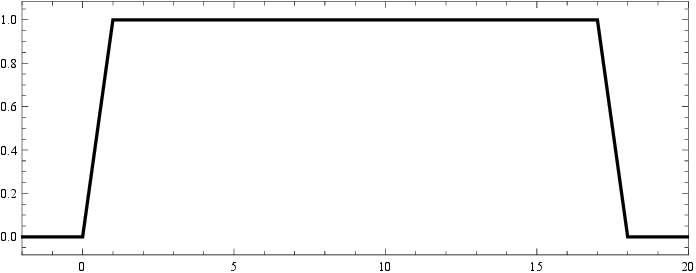}
\end{center}
\caption{A block of the translated autoconvolution $t\mapsto (f_2\ast f_2)(2J_8+t)$. Its plateau has height $1$ and length $16$. Created in Mathematica \cite{Mathematica}.}
\label{fig:m2-convolution}
\end{figure}

At $2J_8=(2n-16)D+(2n+16)D^2$, the even-position digit determines $s=16$. There is exactly one ordered choice with these higher digits in $W_{16}(z)$, since every factor of either type again uses its first position. As $K=16$, \eqref{eq:runs} gives $B_{16}(z)=1+z+\cdots+z^{16}$.
Thus the coefficients of $F(z)^2$ at $2J_8,\ldots,2J_8+16$ are all $1$ and the two immediately adjacent coefficients are zero. On $0\le t\le18$, linear interpolation gives
\[ (f_2\ast f_2)(2J_8+t)=
\begin{cases}
t, & 0\le t\le1,\\
1, & 1<t\le17,\\
18-t, & 17<t\le18 .
\end{cases} \]
Figure~\ref{fig:m2-convolution} shows this function, together with its zero values on the adjacent intervals. The blocks for which all chosen digit positions are distinct have the same profile at level $16$, with every height multiplied by $A_{16}$.

For \underline{$m=1000$} the right-hand side of \eqref{eq:step_est}
is
\[ \left(1-\frac{1}{8\cdot10^{24}}\right) \frac{999}{1000}\left(\frac{999999}{1000000}\right)^2 >0.99. \]
Even the term $U(z)^nV(z)^n$ of $F(z)$ leads to more than $10^{10^{26}}$ steps in the construction of $f_{1000}$. In comparison, the examples of Boyer and Li \cite{BL}, Jaech and Joseph \cite{JJ}, and Georgiev et al.\@ \cite{GGTW} in Table~\ref{tab:records} have 575, 2399, and 50000 steps, respectively.
It is now clear that examples like $f_{1000}$ were unlikely to have been found by numerical searches.


\section*{Declaration of AI usage}
The proof of Theorem~\ref{thm:main} was found in a series of chats with OpenAI's ChatGPT-5.5 Pro, ChatGPT-5.6 Sol, and ChatGPT-6 Astra between June and November 2026. The chats were essentially autonomous, with access to the internet and to the authors' human-written files, including several sketches of the construction used in \cite[\S7.1]{KR26}.
In particular, ChatGPT suggested the polynomial reformulation, the definition of the polynomials $P_i$, the construction of $F$, and the choices of parameters. After ChatGPT-6 Astra arranged the available arguments into a valid proof, the authors simplified it and wrote up the manuscript.
Furthermore, ChatGPT-6 Astra was used for literature searches, copyediting, and proofreading, while the authors assume full responsibility for the content. Finally, ChatGPT-6 Astra was also used to produce the formal verification of Theorem~\ref{thm:main}, which is based on the proof presented in this paper.


\section*{Acknowledgments and funding}
V.\,K. would like to thank Zane Li for a useful discussion on the problem in February 2026.

The paper was supported in part by the Croatian Science Foundation under the project number HRZZ-IP-2022-10-5116 (FANAP).
T.\,C. was also supported in part by the European Union - NextGenerationEU, project: IP-UNIST-44 (ITPEM).
V.\,K. was also supported in part by the European Union - NextGenerationEU through the National Recovery and Resilience Plan 2021--2026, via an institutional grant of the University of Zagreb Faculty of Science, IK IA 1.1.3, \emph{Impact4Math}.


\bibliographystyle{plainurl}
\bibliography{autoconvolutions}

\end{document}